\documentclass[journal,twoside,web]{ieeecolor}
\usepackage{generic}
\usepackage{cite}
\usepackage{amsmath,amssymb,amsfonts}
\usepackage{graphicx}
\usepackage{textcomp}
\usepackage{multicol}
\usepackage{balance}

\usepackage{xcolor}  
\usepackage{cite}
\usepackage{orcidlink}
\usepackage{physics}

\usepackage{tikz}
\usepackage{mathtools}
\usetikzlibrary{arrows,spy,decorations.pathmorphing,decorations.footprints,fadings,calc,trees,mindmap,shadows,decorations.text,patterns,positioning,shapes,matrix,fit,shapes,snakes,3d}
\usetikzlibrary{arrows.meta}

\usepackage{mathtools}
\usepackage{bm}
\usepackage{bbm}

\usepackage{graphicx}
\usepackage{amssymb}

\usepackage{algorithm}
\usepackage{algpseudocode}

\newcommand{\RD}{\color{red}}

\newcommand{\BG}{\mathbb{G}} 
\newcommand{\R}{\mathbb{R}} 
\newcommand{\Bfzr}{\mathbf{0}} 

\newcommand{\CT}{\mathcal{T}} %
\newcommand{\CO}{\mathcal{O}} %
\newcommand{\CG}{\mathcal{G}} %
\newcommand{\CE}{\mathcal{E}} %
\newcommand{\CV}{\mathcal{V}} %

\newcommand{\sign}{{\rm sgn}}

\newcommand{\ctr}{{\rm ctr}}
\newcommand{\dir}{{\rm dir}}

\definecolor{ORG}{HTML}{fe9929}

\newcommand{\blk}{{\rm blkdiag}} 

\newcommand\scalemath[2]{\scalebox{#1}{\mbox{\ensuremath{\displaystyle #2}}}}

\usepackage{amsthm}

\newtheorem{theorem}{Theorem}

\newtheorem{lemma}{Lemma}
\newtheorem{corollary}{Corollary}

\theoremstyle{definition}
\newtheorem{remark}{Remark}

\newtheorem{example}{Example}
\newtheorem{definition}{Definition}
\newtheorem{assumption}{Assumption}

\newenvironment{pfof}[1]{\vspace{1ex}\noindent{\textit{Proof of
		#1:}}\hspace{0.5em}} {\hfill\qed\vspace{1ex}}

\newcommand{\dri}{{\rm drv}}
\newcommand{\bid}{{\rm bid}}
\newcommand{\uni}{{\rm uni}}

\newcommand{\cre}{{\rm cre}}

\newcommand{\inc}{{\rm inc}}
\newcommand{\dec}{{\rm dec}}
\newcommand{\pat}{{\rm pat}}
\newcommand{\cyc}{{\rm cyc}}

\newcommand{\rem}{{\rm rmv}}

\graphicspath{ {./img/} }
\usepackage{hyperref} 
\hypersetup{colorlinks=true,citecolor=[RGB]{5,112,176},linkcolor=[RGB]{5,112,176},urlcolor=[RGB]{5,112,176}}

\def\BibTeX{{\rm B\kern-.05em{\sc i\kern-.025em b}\kern-.08em
    T\kern-.1667em\lower.7ex\hbox{E}\kern-.125emX}}
\begin{document}
\title{Vibrational Control of Complex Networks}
\author{Yuzhen Qin \orcidlink{0000-0003-1851-1370}, Fabio Pasqualetti \orcidlink{0000-0002-8457-8656}, Danielle S. Bassett \orcidlink{0000-0002-6183-4493},  and Marcel van Gerven \orcidlink{0000-0002-2206-9098}
\thanks{Y. Qin and M. van Gerven are with the Department of Machine Learning and Neural Computing, Donders Institute for Brain, Cognition and Behaviour, Radboud University, Nijmegen, the Netherlands (\{yuzhen.qin, marcel.vangerven\}@donders.ru.nl). F. Pasqualetti is with the Department of Electrical Engineering and Computer Science at the University of California, Irvine, CA, 92697, USA, (fabiopas@uci.edu). D. S. Bassett is with the Department of Bioengineering, the Department of Electrical \&
Systems Engineering, the Department of Physics \& Astronomy, the	Department of Psychiatry, and the Department of Neurology, University of Pennsylvania, and The Santa Fe Institute (dsb@seas.upenn.edu). Y. Qin, F. Pasqualetti and D. S. Bassett  were also with  Aligning Science Across Parkinson’s (ASAP) Collaborative Research Network, Chevy Chase, MD 20815. This publication is part of the project Dutch Brain Interface Initiative (DBI2) with project number 024.005.022 of the research programme Gravitation, which is financed by the Dutch Ministry of Education, Culture and Science (OCW) via the Dutch Research Council (NWO). This research was also funded in part by Aligning Science Across Parkinson’s ASAP-020616 through the Michael J. Fox Foundation for Parkinson’s Research (MJFF) and NSF-CMMI-2308639. 
}
\thanks{The manuscript was prepared using LaTex on Overleaf, and the source code can be downloaded at: \hyperref[a]{https://arxiv.org/abs/2408.08263}. All simulations were performed in MATLAB, and the corresponding code is available at this repository \hyperref[b]{https://doi.org/10.5281/zenodo.17169763}.}
}

\maketitle

\begin{abstract}
The stability of complex networks, from power grids to biological systems, is crucial for their proper functioning. It is thus important to control such systems to maintain or restore their stability. Traditional approaches rely on real-time state measurements for feedback control, but this can be challenging in many real-world systems, such as the brain, due to their complex and dynamic nature. This paper utilizes vibrational control, an open-loop strategy, to regulate network stability. Unlike conventional methods targeting network nodes, our approach focuses on manipulating network edges through vibrational inputs. We establish sufficient graph-theoretic conditions for vibration-induced functional modifications of network edges and stabilization of network systems as a whole. Additionally, we provide methods for designing effective vibrational control inputs and validate our theoretical findings through numerical simulations.
\end{abstract}


\section{Introduction}
\label{sec:introduction}
\IEEEPARstart{M}{any} natural and technological systems comprise interacting dynamical units and are often modeled as complex networks.  Effective network systems require stability, expressed as a fixed point, a limit cycle, or a manifold. Maintaining stability is critical for the effective operation of such systems, as instability can have severe consequences, such as catastrophic blackouts in power grids~\cite{JWB-KB-MDB-AQG-EG-VR-JDR-SS-CAS-MEW:21a} and various brain disorders~\cite{JP-DCM-JJGR-SCA:2013}. Consequently, developing control strategies to preserve or restore desired system dynamics is essential.

Most existing control strategies for network systems rely on feedback mechanisms that require real-time access to system states, typically obtained through direct measurement or observers. However, accurately measuring or observing internal states often proves challenging in real-world scenarios. This limitation is particularly pronounced in complex biological systems like the brain, where the dynamic and intricate nature hinders precise state estimation. While techniques such as EEG and MEG provide valuable neural information, their limited temporal and spatial resolution restrict their ability to accurately measure brain activity, especially in real time.

Vibrational control is an open-loop technique that regulates systems without measuring their states. By employing pre-designed, high-frequency dithers, this method has been found to successfully stabilize a wide range of natural and engineered systems, including mechanical systems, chemical reactors, under-actuated robots (see~\cite{REB-JB-SMM:86b,BS-BTZ:97,CX-TY-MI:2018,CK_BB:2024,ZJ-FE:2024} and the references therein), and even insect flight~\cite{THE-KM-HTL-GJSM:2020}. Our previous work suggests that vibrational control might underlie the mechanism of deep brain stimulation, a neurosurgical therapy for brain disorders, due to the shared use of high-frequency stimuli~\cite{YQ-DSB-FP:22a,QY-NA-BDS-PF:2023}. Given the challenges of real-time state acquisition in complex network systems like the brain, vibrational control emerges as a promising approach. This paper initiates an investigation into vibrational control within the framework of linear networks.

\textbf{Related work.} 
Control of complex networks has emerged as a major research area across scientific disciplines. Significant efforts have been directed towards understanding how network structure influences controllability~\cite{YYL-JJS-ALB:11, FP-SZ-FB:13q,JJ-HJW-HLT-KMC:19,GB-KO-AS-KM:21,AW-SM-YY-KX:22} and stabilizability of these systems~\cite{SP-SK-PAA:15,LJ-CX-PS-PGJ-PVM:19}. Various control strategies have emerged, with most focusing on applying feedback control directly to individual nodes within the network~\cite{GLV-FM-LV:2018,GK_RJ_ST:2021TCNS,ZS_ZJ_LJH_LJA:2023,ZLH_WS_LCJ:2022}.  Other approaches aim to directly modify network connectivity to achieve desired behavior~\cite{GLV_FM_et:2020,MS_VTL:2022}. This paper delves into vibrational control as a novel approach for complex network systems. In contrast to existing methods, vibrational control operates in an open-loop manner, introducing control signals to \textit{edges} rather than nodes. This approach draws inspiration from observations in deep brain stimulation, where electrical stimuli primarily affect dendrites and axons near the electrodes, rather than cell bodies~\cite{HMH-BV-RRC:09}. A paper on bilinear systems shares some similarities with our work~\cite{Gharesifard_auto:2017}. It focuses on extending the system's graph using vibration-like inputs, such that the extended graph admits a Hurwitz weighted adjacency matrix when appropriately weighted. By contrast, our paper aims to directly stabilize given network systems.

\textbf{Contribution.} The contributions are fourfold. First, we introduce the ``functioning network'', a novel concept that captures the average behavior of a network system under vibrational control. This concept elucidates the working mechanism of vibrational control in network systems: indirectly modifying network connections. Second, we define four new edge properties, \textit{increasability, decreasability, removability, and creatability,} characterizing how individual edges can be manipulated through vibrational control. We further derive graph-theoretic conditions for these properties.  Notably, by introducing the concept of joint vibrations with identical frequencies, our conditions are much less restrictive than existing ones on general linear systems. Third, we present conditions for simultaneously modifying multiple edges without unintended side effects, paving the way for designing sophisticated control strategies. Finally, we establish sufficient conditions for vibrational stabilization of entire network systems and provide design principles for constructing effective vibrational inputs. Finally, numerical experiments are  conducted to demonstrate our theoretical findings.   

This paper substantially extends our previous work~\cite{QY-NA-BDS-PF:2023} on characterizing the vibrational modifiability of network edges, particularly by leveraging same-frequency vibrations. Consequently, our results on stabilizing network systems are significantly more comprehensive.

\textbf{Notation}: Given any matrix $A\in\R^{n\times n}$, one can associate it with a weighted directed graph (digraph), denoted as $\CG(A):=(\CV,\CE,A)$. Here, $\CE=\{1,2,\dots,n\}$ is the node set, and $\CE\subseteq \CV\times \CV$ is the edge set. A direct edge from $i$ to $j$, denoted as $(i,j)$, satisfies $(i,j)\in\CE$ if and only if $a_{ji}\neq 0$;  $a_{ij}$ is its weight. The sign graph associated with $A$ is denoted as $\CG_{\sign}(A):=(\CV,\CE,\sign(A))$. Given two sign graphs $\CG_1=(\CE_1, \CV, S_1)$ and $\CG_2=(\CE_2, \CV, S_2)$, we denote $\CG_1\subseteq \CG_2$ if $\CE_1\subseteq \CE_2$ and each edge in $\CE_1$ has the same sign as their counterpart in $\CE_2$. In addition, we let $\BG(A)$ be the unweighted digraph associated with $A$. A continuous function $\alpha: [0,\infty)\to [0,\infty)$ is said to belong to function class $\mathcal {K}$ if it is strictly increasing and satisfies $\alpha(0)=0$.

\section{Problem Formulation}

Consider a linear network system  governed by 
\begin{equation}\label{main}
	\dot x_i(t) = a_{ii} x_i(t)+\sum_{j=1, j\neq i}^{n}a_{ij} x_j(t), 
\end{equation}
where $x_i\in\R, i=1,2,\dots,n$, is the state of the $i$th subsystem,  $a_{ii}\in\R$ represents its intrinsic dynamics, and $a_{ij}$ describes the interconnection from subsystems $j$ to $i$. 
Defining $x=[x_1,x_2,\dots,x_n]^\top\in \R^n$ and $A:=[a_{ij}]_{n\times n}$, the  system~\eqref{main} can be rewritten into the  compact form:
\begin{equation}\label{main:compact}
	\dot x = A x.
\end{equation}
It can be observed that the network of this system is described by the weighted digraph $\CG :=\CG(A)$.

In this paper, we assume that the system~\eqref{main:compact} is unstable (i.e., $A$ is not Hurwitz), and our objective is to design controllers to stabilize this system. In particular, we are interested in stabilizing the system using vibrational control. Next, we introduce the concept of vibrational control.

\subsection{Vibrational Control}

Given a general linear system of the form~\eqref{main:compact}, with $A$ not necessarily formed by a network, consider a control matrix $W(t)=[w_{ij}(t)]_{n\times n}$ that influences the system parameters in $A$, resulting in the  controlled system below:
\begin{equation}\label{controlled_net_compact}
	\dot x = \big(A+W(t) \big) x. 
\end{equation} 
The input matrix $W(t)$ is  often chosen to have the form of 
\begin{align}\label{vib:form}
	w_{ij}(t)= \sum_{\ell=1}^{\infty} \alpha_{ij}^{(\ell)} \sin\left(\ell \beta_{ij}t+\phi_{ij}^{(\ell)}\right),
\end{align}
where $ \alpha_{ij}^{(\ell)}, \phi_{ij}^{(\ell)}\in\R$, and $\beta_{ij}>0$ for any $i, j$, and $\ell$. Each $w_{ij}(t)$ is (almost) periodic and zero-mean~\cite{REB-JB-SMM:86b}, satisfying
$
 \lim\limits_{T\to \infty}\frac{1}{T}\int_{t=0}^{T}w_{ij}(t) \dd t=0.
$
A wide range of periodic signals, including square, triangular, and sinusoidal waves, can be represented or approximated by~\eqref{vib:form}.
Typically, each $w_{ij}(t)$ is high-frequency, introducing vibrations to the parameter $a_{ij}$. Therefore, control input $W(t)$ of the form~\eqref{vib:form} is termed \textit{vibrational control}. By carefully configuring these vibrations, unstable systems can be stabilized without requiring state measurements~\cite{SMM:80,BS-BTZ:97,CX-TY-MI:2018}. A recent intriguing study reveals that high frequencies are not necessary; however, as the frequency decreases, the range of stabilizing parameters narrows~\cite{CK_BB_2025}. In this paper, we still focus on high-frequency vibrations.

\subsection{Vibrational Control in Network Systems}
For general linear systems, vibrations can be introduced to any $a_{ij}$ in the system matrix $A$. However, for network systems, we assume that vibrations can only be introduced to the nonzero entries in $A$, corresponding to the existing edges in the network. This assumption is motivated by the possible connections between vibrational control and brain stimulation, as discussed in Section~\ref{sec:introduction}. The direct influence of electrical stimuli delivered by stimulation electrodes is naturally limited to existing neural connections. Formally, we have the following constraint on the vibrational control matrix $W(t)$:
\begin{align}\label{constrain}
	&w_{ij}(t)=0, \forall t\ge 0,& \text{ whenever } a_{ij}=0.
\end{align}
In other words, the non-zero pattern of the vibrational control matrix $W(t)$ must be constrained by that of the matrix $A$. Notably, this constraint introduces additional technical challenges compared to existing work without such a restriction~\cite{SMM:80}, as control placements in our case are less flexible.
For clarity, throughout this paper, ``vibrational control'' refers exclusively to control inputs satisfying Eqs.~\eqref{vib:form} and~\eqref{constrain}.

Given a vibrational control described by $W(t)$, we say the network system~\eqref{main:compact} is vibrationally \textit{stabilized} if the controlled system~\eqref{controlled_net_compact} is asymptotically stable.

\begin{definition}[Vibrational stabilizability]
	The network system~\eqref{main:compact} is said to be \emph{vibrationally stabilizable} if there exists a vibrational control that stabilizes the system~\eqref{main:compact}. 
\end{definition}

\section{Averaged Systems and Functioning Networks}\label{functioning}
To explicitly highlight the high-frequency nature of the vibrational inputs, without loss of generality, we rewrite the system~\eqref{controlled_net_compact} as follows:
\begin{equation}\label{system:vib}
	\dot x = \left(A+\frac{1}{\varepsilon}V\left(\frac{t}{\varepsilon}\right) \right) x, 
\end{equation}
where $V(t/\varepsilon)/\varepsilon=W(t)$, and $\varepsilon>0$ determines the frequencies of the vibrations. Smaller values of $\varepsilon$ correspond to higher vibration frequencies. Consequently, it becomes to design $V(t)$ and the parameter $\varepsilon$ to stabilize the original system~\eqref{main:compact}.

To assess the stability of a network system under vibrational control, the controlled system~\eqref{system:vib} needs to be analyzed. A common approach is to apply averaging techniques (e.g.,~\cite{SMM:80,QY_KY_ABDO_CM:2021,BS-BTZ:97}). To do that, we first introduce a new timescale $s={t}/{\varepsilon}$.  This transformation leads to 
\begin{equation}\label{system:timescale}
	\frac{\dd x}{\dd s} = (\varepsilon A + V(s))x. 
\end{equation}

Since $V(s)$ has a zero mean, the conventional averaging method (see~\cite[Chap. 10]{HKK:02-bis} and~\cite{SJA-FMJ:07}) cannot be directly applied to the system~\eqref{system:timescale} as it would simply eliminate the $V(s)$ term, reverting the system to the uncontrolled form ${d x}/{d s} = \varepsilon A x$.
To avoid this, we introduce a coordinate transformation to the system~\eqref{system:timescale} before applying the averaging technique. Specifically, we define an auxiliary system:
\begin{equation}\label{auxiliary}
	\dot{ \hat x} = V(s) \hat x.
\end{equation}
Let $\Phi(s,s_0)$ be the state transition matrix of~\eqref{auxiliary} satisfying $\hat x(s)=\Phi(s,s_0)x(s_0)$. For the state transition matrix, there always exists a fundamental matrix $\Psi(t)$ such that $\Phi(s,s_0)= \Psi(s) \Psi^{-1}(s_0)$ (see~\cite[Ch. 2.3]{AP-MAN:97} for the definition of fundamental matrix). 
Introducing a change of variables to the system~\eqref{system:timescale} by letting $z=\Psi^{-1}(s) x$, one can derive that
\begin{equation}\label{model:changed_coord}
	\frac{\dd z}{\dd s}=\varepsilon \Psi^{-1}(s) A  \Psi(s) z. 
\end{equation}
Since $V(t)$ is almost periodic and zero-mean, $\|\Psi(s)\|$ and $\|\Psi^{-1}(s)\|$ are always bounded. Then, $x=0$ of the system~\eqref{system:vib} is asymptotically stable if $z=0$ of the system~\eqref{model:changed_coord} is. 

Now, let us consider the averaged system of~\eqref{model:changed_coord}
\begin{equation}\label{averaged}
	\frac{\dd \bar x}{\dd s}= \varepsilon \bar A \bar x, 
\end{equation}
where $\bar A =[\bar a_{ij}]_{n\times n}$ is given by
\begin{equation*}
	\bar A =\lim_{T\to \infty}\frac{1}{T} \int_{s=0}^{T} \Psi^{-1}(s,s_0) A \Psi(s,s_0) d s.
\end{equation*}

The following lemma, which directly follows from Theorem~10.4 in~\cite{HKK:02-bis}, connects the behavior of the averaged system~\eqref{averaged} to that of the  system~\eqref{system:timescale}.
\begin{lemma}\label{lemma:averaged}
	Let $x(t)$ and $\bar x (t)$ be the solutions to the systems~\eqref{system:timescale} and~\eqref{averaged}, respectively, with the initial conditions satisfying $x(0)=\bar x(0)$. Then, the following statements hold:
	\begin{enumerate}
		\item[(i)] There exist a $\mathcal {K}$ function $\kappa$ and $T>0$ such that $\|x(t)-\bar x (t) \|\le \CO(\kappa(\varepsilon))$ for any $t \in [0,T/\varepsilon]$.  
		\item[(ii)] If $\bar x=0$ of the averaged system~\eqref{averaged} is exponentially stable, then there exists $\varepsilon_0>0$ such that, for any $0 < \varepsilon < \varepsilon_0$,  the system~\eqref{system:timescale} is also exponentially stable. 
	\end{enumerate}
\end{lemma}

This lemma states that the solution of the averaged system~\eqref{averaged} provides an $\CO(\kappa(\varepsilon))$ approximation for the solution of the  system~\eqref{system:timescale} over the time interval $[0,T/\varepsilon]$. Furthermore, if the averaged system is stable, so is the original one if the frequency of the vibrational inputs is sufficiently high. 

 Changing back the timescale of~\eqref{averaged} to $t=\varepsilon s$, we arrive at 
\begin{equation}\label{averaged:main}
	\dot {\bar x} = \bar A \bar x.
\end{equation}

\begin{figure}[t]
	\centering
	\begin{tikzpicture}[scale=1]
		\node[inner sep=0pt, scale=0.7] (f1) at (0,0.0) {$
			\dot x =A x 
			=
			\begin{bmatrix*}[r]
				-1 & 0 & 0  & 0  & 0\\
				0  & 0& 0  & 0  & -0.5 \\
				2  & 2 & 0 & -0.2 & 0 \\ 
				0  & 0 & 0 & 1  & 3 \\
				2 & 0 & -1  & 0  & -2
			\end{bmatrix*}
			x$};
		\node[inner sep=0pt, scale=1.05] (f1) at (4.5,0.2) {\includegraphics{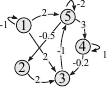}};
		
		\node[inner sep=0pt, scale=0.9] (f1) at (-2,0.6) {(a)};
		
		\node[inner sep=0pt, scale=0.9] (f1) at (-2,-1.3) {(b)};
		\node[inner sep=0pt, scale=0.9] (f1) at (-2,-2) {(c)};
		
		\node[inner sep=0pt, scale=0.7] (f1) at (4.5,-0.8) {Original network $\CG=\CG(A)$};
		
		\node[inner sep=0pt, scale=0.7] (f1) at (-0.5,-1.3) {$
			\dot x =(A+V(t/\varepsilon)/\varepsilon) x $};
		\node[inner sep=0pt, scale=0.7] (f1) at (3,-1.3) {Vibrationally controlled system};
			
		\draw[dashed,line width=0.1mm] (-2.2,-1) -- (6,-1); 
		\draw[-{Stealth[scale=0.8]},line width=0.1mm] (0,-1.48)  -- (0,-1.98); 
		
			\draw[dashed,line width=0.1mm] (-2.2,-1.55) -- (6,-1.55); 
		
		\node[inner sep=0pt, scale=0.7] (f1) at (0.57,-1.76) {averaging};
			
		\node[inner sep=0pt, scale=0.7] (f1) at (-0.2,-2.9) {$
			\dot {\bar x} =\bar A \bar x 
			=
				\begin{bmatrix*}[r]
				-1 & 0 & 0  & 0  & 0\\
				0  & 0& 0  & 0  & -0.5 \\
				{\RD 0}  & 2 & 0 & {\RD 1}& {\RD 2} \\ 
				0  & 0 & 0 & 1  & 3 \\
				{\RD 1.2} & 0 & -1  & 0  & -2
			\end{bmatrix*}
			x$};

		\node[inner sep=0pt, scale=1.05] (f1) at (4.5,-2.5) {\includegraphics{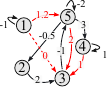}};
		\node[inner sep=0pt, scale=0.7] (f1) at (4.5,-3.5) {Functioning network $\bar \CG=\CG(\bar A)$};
		
	\end{tikzpicture}
	\caption{Schematic of the averaged system and functioning network. (a) The original system and the associated network. (b) The system under vibrational control. (c) The averaged system and the functioning network associated with it. Compared to the original one, some edges have been increased in weight [$(4,3)$], decreased in weight [$(1,5)$],  completely removed [$\RD (1,3)$], and created [$(5,3)$] in the functioning network. }
	\label{func_net}
\end{figure}
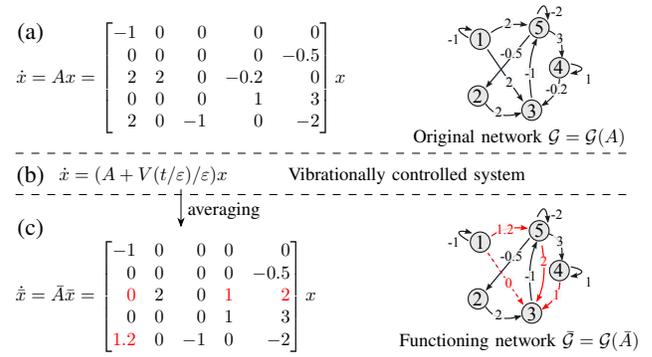

 By comparing systems~\eqref{main:compact} and~\eqref{averaged:main}, we can interpret vibrational control as modifying the system matrix from $A$ to $\bar A$ on average.
This is the underlying mechanism of how vibrational inputs are able to  stabilize a system~\cite{SMM:80,BS-BTZ:97,CX-TY-MI:2018}.

\textbf{Functioning network.} Let us associate the averaged system~\eqref{averaged:main} with a network described by a weighted digraph $\bar \CG:= \CG(\bar A)=(\CV,\bar\CE, 
\bar A)$.  We refer to this network described by $\bar \CG$ as the \textit{functioning network} as it reflects the system's average functionality, at least in the time interval $[0,T/\varepsilon]$. Figure~\ref{func_net} illustrates the relation between the original network $\CG$ and the functioning one $\bar \CG$. 
One can discern that vibrational control can induce changes to the network of the system~\eqref{system:vib} by modifying its connection weights or  altering its network structure.  

Stabilizing an inherently unstable network system~\eqref{main:compact} requires precise implementation of such modifications. First, in Sections~\Ref{controllability} and~\Ref{multiple_contr}, we investigate how vibrational control can reliably induce changes at both the individual edge level and across multiple edges simultaneously. Then, in Section~\Ref{stabilizability}, we explore how to achieve overall network system stabilization by strategically integrating these modifications.

\section{Vibrational Modification of Individual Edges}\label{controllability}
We explore how individual edges of the network $\CG (A)$ associated with the system~\eqref{main:compact} can be vibrationally modified. 

For an existing edge $(j,i) \in \CE$, vibrational control can induce three types of modifications (e.g., see Fig.~\ref{func_net}):
\begin{enumerate}
	\item[(i)]  increase in weight: $\bar a_{ij}>a_{ij}$,
	\item[(ii)] decrease in weight: $\bar a_{ij}<a_{ij}$,
	\item[(iii)] edge removal: $\bar a_{ij}=0$ while  $ a_{ij}\neq 0$.
\end{enumerate}
Additionally, vibrational control can potentially create a new edge that does not exist in the original network (i.e., $\bar a_{ij}\neq 0$ while $a_{ij}=0$). We collectively refer to all these alterations as \textit{functional changes}. Formal definitions are provided below.

\begin{definition}[Vibrational increasability, decreasability, and removability]\label{existing}
	Consider an edge $(j,i)\in \CE$.  We say it is \emph{vibrationally increasable}, (resp., \emph{decreasable}) if for any $\delta> 0$, there exists a vibrational control matrix $V(t)$ such that $\bar a_{ij}=a_{ij}+\delta$ (resp., $\bar a_{ij}=a_{ij}-\delta$). It is said to be \emph{vibrationally removable} if there exists $V(t)$ such that $\bar a_{ij} = 0$. 
\end{definition}

Note that an edge can be simultaneously increasable and decreasable; we refer to such an edge as \textit{vibrationally controllable}\footnote{This definition of controllability is slightly different from the one in~\cite{SMM:80}, where an entry in $A$ is called controllable as long as its value can be functionally changed (increased or decreased), not necessarily in both directions.}. A controllable edge is also removable since it can be removed by simply increasing its weight by $|a_{ij}|$ if $a_{ij}<0$ or decreasing by $a_{ij}$ if $a_{ij}>0$. 

\begin{definition}[Vibrational creatability]\label{no_existing}
	Consider an edge $(j,i)\notin\CE$. We say it is \emph{vibrationally creatable} if for any $\delta \neq  0$, there exists a vibrational control matrix $V(t)$ such that $\bar a_{ij} =\delta$.
\end{definition}

For simplicity, we will henceforth refer to the edges in Definitions~\ref{existing} and \ref{no_existing} as increasable, decreasable, removable, or creatable edges, dropping the ``vibrationally'' prefix.
We proceed by constructing conditions for them, alongside controller design to achieve these modifications. Notably, some edges can be altered through direct vibrational input, while others necessitate the combined effect of multiple vibrations. We delve into these distinct scenarios separately.

\subsection{Modifying Edges via Direct Vibrations} \label{sec:direct_mofi}

\begin{theorem}[Sufficient conditions for direct increasability and decreasability]\label{ctr_bility}
	For the digraph $\CG=(\CV,\CE,A)$ associated with the network system~\eqref{main:compact}, an edge $(j,i)\in \CE$ is 
	
	\begin{enumerate}
		\item[(i)]   increasable if the edge in the opposite direction exists, i.e.,  $(i,j)\in \CE$, and is negatively weighted ($a_{ji}<0$);
		\item [(ii)]  decreasable if the edge in the opposite direction exists and is positively weighted ($a_{ji}>0)$. 
	\end{enumerate}
\end{theorem}

The next lemma presents a method for designing vibrational control inputs to arbitrarily increase or decrease the weight of an edge that meets the above conditions. 
\begin{lemma}[Control design]\label{design:direct}
	Assume that the edge $(j,i)\in \CE$ satisfies condition (i) in Theorem~\ref{ctr_bility} (resp., condition (ii)). Then, for any $\delta>0$, the vibrational control matrix  $V(t)=[v_{k\ell}(t)]_{n \times n}$ in the system~\eqref{system:vib} satisfying the following conditions functionally increases (resp., decreases) the weight of $(j,i)$ by $\delta$:
    \begin{align}\label{design:vibra_form}
		v_{k\ell}(t)	= \begin{cases}
			u_{k\ell}\sin (\beta_{k\ell} t ), \hspace{.2cm}\text{ if } k=i, \ell =j,\\
			0,  \forall t\ge 0, \hspace{0.7cm}\text{ otherwise},
		\end{cases}
	\end{align}
	where the amplitude $ u_{k\ell}$ and frequency (determined by $\beta_{k\ell}$) are arbitrarily chosen such that $\frac{u_{ij}}{\beta_{ij}}=\sqrt{{2 \delta}/{|a_{ji}|}}.$
	\end{lemma}

The analysis of Theorem~\ref{ctr_bility} and Lemma~\ref{design:direct} follows the same steps as that in~\cite{SMM:80}.  Further, we provide a corollary establishing a sufficient condition for vibrational edge removal and outline a corresponding control design approach.

\begin{corollary}\label{coro:rem}
	An edge $(j,i)\in \CE $ is  removable if the edge in the opposite direction exists and has the same sign, i.e., $ \sign(a_{ij})= \sign(a_{ji})$. The vibrational control matrix $V(t)$ functionally removes $(j,i)$ if it satisfies  the condition~\eqref{design:vibra_form} with arbitrarily chosen amplitude and frequency satisfying
		$\frac{u_{ij}}{\beta_{ij}}=\sqrt{{2 a_{ij}}/{a_{ji}}}.$
\end{corollary}

 \begin{figure}[t]
	\centering
	\includegraphics[scale=1.05]{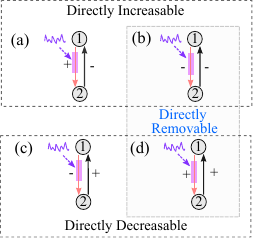}
	\caption{Illustration of directly modifiable edges. The red edge $(1,2)$ can be functionally increased (resp. decreased) by direct vibrational inputs if the opposite edge $(2, 1)$ exists and is negatively (resp. positively) weighted. It is directly removable if the reverse edge has the same sign.}
	\label{dir_removale_types}
\end{figure} 
 
Lemma~\ref{design:direct} and Corollary~\ref{coro:rem} together demonstrate that, to modify edges under the respective conditions, it only requires vibrations applied directly to those edges. Therefore, we refer to these edges  as \textit{directly increasable}, \textit{decreasable}, or \textit{removable edges}. Collectively, we call them \textit{directly modifiable}. Fig.~\ref{dir_removale_types} provides an illustration of these edges and the constructed sufficient conditions for them. Notably, directly removable edges always come in pairs. This is because if edge $(j,i)$ fulfills the conditions of Corollary~\ref{coro:rem}, the edge in the opposite direction, $(i,j)$, naturally satisfies them as well.
  
 \begin{remark}
 	As in Lemma~\ref{design:direct} and Corollary~\ref{coro:rem}, this paper focuses on designing sinusoidal vibrations  that naturally satisfy~\eqref{vib:form} for inducing desired functional changes. Other signal waveforms, such as square or triangular waves as discussed in~\cite{Gharesifard_auto:2017}, could potentially be employed for the same purposes.
 \end{remark}

\subsection{Modifying Edges via Joint Vibrations} \label{joint_control}
In this subsection, we present conditions on which edges can be modified by combining multiple vibrational inputs. Before that, let us provide a useful graph-theoretic definition.

\begin{definition}[Directed trail and path~\cite{BJA-MUSR:1976}]
	Given a digraph $\CG=(\CV,\CE,A)$, a directed \textit{trail}  is a sequence of edges $\{e_1, e_2, \dots, e_k\}$ such that: (i)  each $e_i =(v_i,v_{i+1})\in\CE$ for $i=1, 2, \dots, k$, and (ii) all edges are distinct. A trail with all distinct nodes is called a \textit{path}. If $v_1=i$ and $v_{k+1}=j$, we say it is a trail from $i$ to $j$, denoted as $\CT_{ij}$. 
\end{definition}

\begin{theorem}[Sufficient condition for joint controllability]\label{ctr:indirect}
	For the graph $\CG=(\CV,\CE,A)$ associated with the network system~\eqref{main:compact}, an edge $(j,i)\in \CE$ is vibrationally controllable if there are two nodes $p,q\in \CV$ such that $\CT_{ji}=\{(j,p), (p,q), (q,i)\}$ is a directed trail that contains no self-loop.
\end{theorem}

\begin{figure}[t]
	\centering
	\includegraphics[scale=1.05]{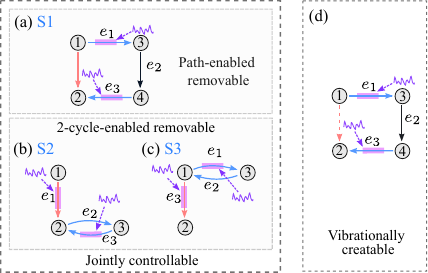}
	\caption{Illustration of jointly controllable and creatable edges. (a)-(c) Jointly controllable edges:	each corresponds to the situations~\eqref{sit_1}-\eqref{sit_3}, respectively. Injecting appropriate same-frequency vibrations into the driver edges (i.e., $e_1$ and $e_3$) can induce an arbitrary change to $(1,2)$. A controllable edge is naturally removable; removable edges are further classified as path-enabled (a) or 2-cycle-enabled (b-c) based on the connection structure.  (d) Jointly creatable edges: joint vibrations can also create an arbitrarily-weighted edge. }
	\label{combi_contr_types}
\end{figure} 

There are only three possible situations that satisfy the conditions in this theorem (see Fig.~\ref{combi_contr_types} (a)-(c) for an illustration):
		\begin{align}
		&\text{(i)}~ \text{The nodes } j, p, q, i \text{ are all distinct},\tag{S1} \label{sit_1}\\
		&\text{(ii)}~ p=i,~\text{where}~\CT_{ij}~ \text{becomes} \{(j,i), (i,q), (q,i)\} \tag{S2} \label{sit_2},\\
		&\text{(iii)}~ q=j,~\text{where}~\CT_{ij}~ \text{becomes} \{(j,p), (p,j), (j,i)\} \tag{S3} \label{sit_3}.
	\end{align}

Note that edges satisfying the condition in Theorem~\ref{ctr:indirect} is simultaneously increasable and decreasable. The lemma below outlines a method for designing joint inputs to achieve desired changes, thereby proving Theorem~\ref{ctr:indirect}. 

\begin{lemma}\label{indirect:design}
	Consider an edge $(j,i)\in \CE$ in $\CG$ satisfying the condition in Theorem~\ref{ctr:indirect}. Consider the vibrational control matrix  $V(t)=[v_{k\ell}(t)]_{n \times n}$ in the system~\eqref{system:vib} with 
\begin{align}\label{theo:amp}
		&v_{pj}(t)	=u_{pj}\sin (\beta t ), &v_{iq}(t)	= u_{iq}\sin (\beta t ),
	\end{align}
 and $v_{k\ell}(t)=0, \forall t \ge 0,$ for any $k,\ell$ that satisfy $k\neq i$ or  $\ell \neq j$. Then, for any $\delta \in \R$, the weight of $(j,i)$ is changed to $\bar a_{ij} =a_{ij}+\delta$ if $u_{pj}, u_{iq}$, and $\beta$ are chosen such that	
	\begin{align*}
		&\text{(i)}~ \delta= -\frac{1}{2\beta^2}  u_{iq} a_{qp} u_{pj}~\text{in Situation~\eqref{sit_1}}, \\
		&\text{(ii)}~ \delta=- \frac{1}{2\beta^2} \big( u_{iq} a_{qi} u_{ij}  + u_{ij}^2 a_{ji} \big)~\text{in Situation~\eqref{sit_2}}, \\
		&\text{(iii)}~ \delta=- \frac{1}{2\beta^2} \big(u_{ij} a_{jp} u_{pj}   +  u_{ij}^2 a_{ji} \big)~\text{in Situation~\eqref{sit_3}}.
	\end{align*}
\end{lemma}

Recall that a controllable edge is also removable, meaning that edges under the conditions~\eqref{sit_1}-\eqref{sit_3} are removable. We refer to such edges as \textit{jointly removable} edges to distinguish from directly removable ones defined in Section~\ref{sec:direct_mofi}. In Situation~\eqref{sit_1}, $(j,i)$ is removable due to the existence of a directed path from $j$ to $i$ since $i, j, p, q$ are all distinct.  Therefore, we refer to $(j,i)$ as \textit{path-enabled removable edge}. In Situations~\eqref{sit_2} and~\eqref{sit_3}, $(j, i)$ is removable due to its connection to a cycle with two nodes. We refer to such an edge as~\textit{2-cycle-enabled removable edge}. 
In all the  situations~\eqref{sit_1}-\eqref{sit_3}, the edge $(j,i)$ can be functionally modified by simultaneously introducing vibrations to the edges $(j,p)$ and $(q,i)$. We refer to these two edges as the \textit{driver edges} of $(j,i)$. This is distinct from the case for a directly modifiable edge, where the driver edge is simply itself since a direct vibration can functionally change it. Fig.~\ref{combi_contr_types} presents a more intuitive illustration of these concepts.

The next corollary that  follows from Theorem~\ref{ctr:indirect} and Lemma~\ref{indirect:design} presents a condition for vibrational creatability.
\begin{corollary}[Creatability]\label{coro:creatable}
	The edge $(j,i)\notin \CE$  is vibrationally creatable if there exists two nodes $p$ and $q$ such that $\{(j,p), (p,q), (q,i)\}$  is a directed path. In addition, for any $\delta \in \R$, the edge $(j,i)$  is created with weight $\bar a_{ij} =\delta$ if $u_{pj}, u_{iq}$, and $\beta$ in~\eqref{theo:amp} are chosen such that	
	$\delta= -\frac{1}{2\beta^2}  u_{iq} a_{qp} u_{pj}.$
\end{corollary}

The joint vibrational control in this subsection necessitates \textit{identical frequencies} for vibrations applied to both edges. This approach distinguishes from existing methods that typically employ incommensurable frequencies\footnote{Given two non-zero real numbers $a,b\in\R$, they are said incommensurable if their ratio $a/b$ is not a rational number.} to isolate effects of vibrations introduced to different locations of a system~\cite{SMM:80}. Interestingly, vibrations with identical frequencies can control edges that prove challenging to control via direct vibrations.

	In cases~\eqref{sit_1}-\eqref{sit_3}, edge $(j,i)$ can potentially be directly modified if $a_{ji}\neq 0$. 	
	While direct vibration can alter this edge, the modification is constrained to either increasing or decreasing its weight, contingent on the sign of the opposite edge. In contrast, applying joint, same-frequency vibrational inputs offer greater flexibility, enabling the edge's weight to be increased, decreased, and removed,  irrespective of the opposite edge's sign. In addition, same-frequency joint vibrations also enable the creation of new edges. Similar ideas are considered in~\cite{Gharesifard_auto:2017} to functionally extend the network in bilinear systems. 
	
Despite its advantages, joint vibrational control can introduce unintended modifications due to network effects.   The following example illustrates this phenomenon.

\begin{figure}[t]
	\centering
	\begin{tikzpicture}[scale=1]

		\node[inner sep=0pt, scale=0.7]  at (-4.5,3.9) {(a) };
		\node[inner sep=0pt, scale=0.7]  at (-3,3.9) {Original system};
		\node[inner sep=0pt, scale=0.7]  at (-4.5,1.1) {(b)};
		\node[inner sep=0pt, scale=0.7]  at (-2.1,1.1) {Controlled system with a direct vibration};
		\node[inner sep=0pt, scale=0.7]  at (-4.5,-1.5) {(c)};
		\node[inner sep=0pt, scale=0.7]  at (-2.25,-1.5) {Controlled system with joint vibrations};
		\node[inner sep=0pt, scale=0.65] at (2,0) {\includegraphics{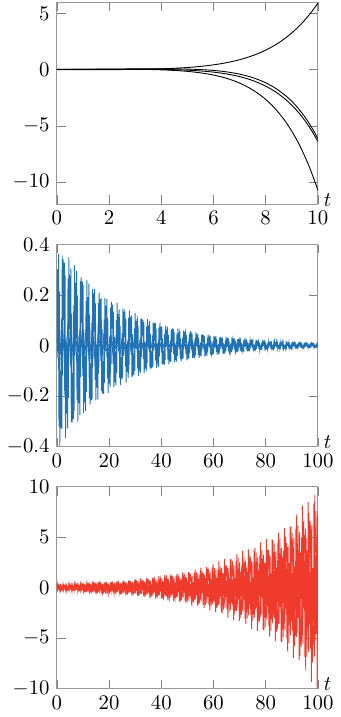}};
		
		\node[inner sep=0pt, scale=0.7]  at (-1,-.95) {Functioning};
		
		\node[inner sep=0pt, scale=0.7]  at (-3.5,-.95) {Driver edge};
		
		\node[inner sep=0pt, scale=0.7]  at (-3.5,-3.9) {Driver edges};
		\node[inner sep=0pt, scale=0.7]  at (-1,-3.9) {Functioning};
		
		\node[inner sep=0pt, scale=0.7]  at (-2.3,2.7) {$
			A
			=
			\begin{bmatrix*}[r]
				0.1 &0 &0 &-1\\
				1 &-1 &0 &0 \\
				0&1 &-0.3&0 \\
				-1 &0 &1 &-0.2
			\end{bmatrix*} 
			$};
		\node[inner sep=0pt, scale=1.05] at (-3.6,0) {\includegraphics{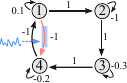}};
		\node[inner sep=0pt, scale=1.05] at (-1.1,0) {\includegraphics{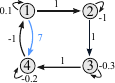}};
		
		\node[inner sep=0pt, scale=1.05] at (-3.6,-2.6) {\includegraphics{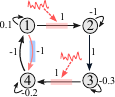}};
		\node[inner sep=0pt, scale=1.05] at (-1.1,-2.75) {\includegraphics{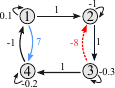}};		
	\end{tikzpicture}
	\caption{Joint vibrational control can cause unintended functional changes. (a) The original system is unstable. (b) A direct vibrational control stabilizes the system. (c) A joint vibrational control fails to stabilize the system due to an unintended creation of the edge $(3,2)$.}
	\label{single_edge_stabilization}
\end{figure}
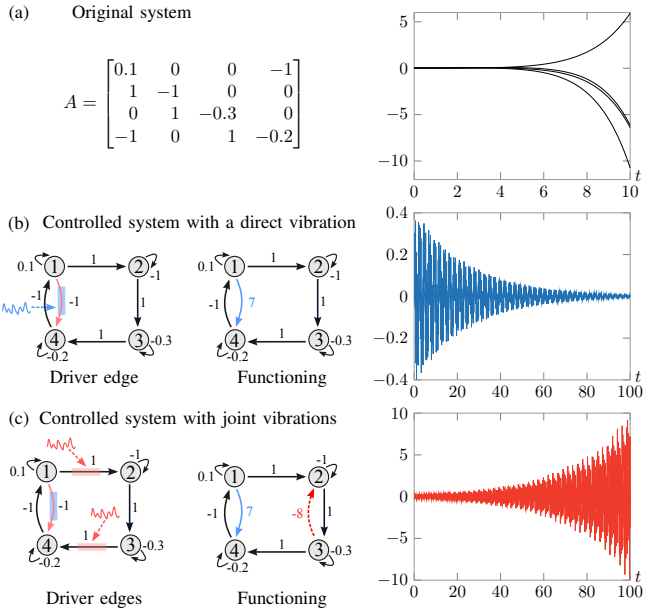

\begin{example}
	Consider an unstable network system with $A$ given in Fig.~\ref{single_edge_stabilization}. It can be observed that the edge $(1,4)$ is both directly increasable and jointly controllable. By calculation, functionally increasing the weight of $(1,4)$ from $-1$ to $7$ leads to a stable averaged system. To realize this change, there are two approaches for vibrational control: (i) introduce a direct vibration to the edge $(1,4)$ or (ii) jointly inject same-frequency vibrations to the edges $(1,2)$ and $(3,4)$.   
	
	Leveraging Lemma~\ref{design:direct}, we let the direct vibrational input be  $v_{41}(t)=4 \sin(\beta t )$. As shown in Fig.~\ref{single_edge_stabilization}-(b), the system is effectively stabilized (where $\varepsilon=0.04$).	
	
	For joint vibrational control, utilizing Lemma~\ref{indirect:design}, we let the joint vibrations injected into $(1,2)$ and $(3,4)$ be 
	$
		v_{21}(t) = 4 \sin(\beta t),  v_{43}(t) = -4 \sin(\beta t)
	$
 with the same $\varepsilon$.
	While this joint vibrational control achieves the desired functional change on edge $(1,4)$, it fails to stabilize the system as shown in Fig.~\ref{single_edge_stabilization}-(c). This is because the simultaneous vibrations applied to the two driver edges create an unintended side effect: the creation of a new edge $(3,2)$. 
\end{example}

Stabilizing a system can sometimes be achieved by controlling a single edge, as illustrated in Fig.~\ref{single_edge_stabilization}-(b). However, in most cases, controlling multiple edges simultaneously is necessary. It is important to note that applying multiple vibrations even with different frequencies to the network simultaneously can induce unintended modifications due to  network effects, similar to the issue shown in Fig.~\ref{single_edge_stabilization}-(c). Therefore, careful selection of driver edges and precise design of control inputs are crucial for achieving desired functional changes.

\section{Simultaneous Modification of Multiple Edges}\label{multiple_contr}

This section focuses on simultaneously modifying multiple edges. To facilitate this analysis, we define several edge sets within the digraph $\CG$ associated with the system~\eqref{main:compact}. These sets categorize edges based on their modifiability through vibrational control, i.e., functional increase, decrease, bidirectional control (increase and decrease), or creation.

\begin{definition}
	Consider the digraph $\CG=(\CV,\CE,A)$ associated with the system~\eqref{main:compact}. We define the following edge sets:
	
	\begin{enumerate}
		\item[(i)] The directly increasable edge set $\CE_\inc \subseteq \CE$ is the set of edges  satisfying the condition in Theorem~\ref{ctr_bility}-(i);
		\item[(ii)] The directly decreasable edge set $\CE_\dec \subseteq \CE$ is the set of edges  satisfying the condition in Theorem~\ref{ctr_bility}-(ii);	
		\item[(iii)]  The controllable edge set $\CE_\ctr \subseteq \CE$ is the set of edges satisfying the conditions in Theorem~\ref{ctr:indirect};
		\item[(iv)]  The creatable edge set $\CE_\cre\nsubseteq \CE$ with $\CE_\cre\subseteq \CV \times \CV$ is the set of edges  satisfying the condition in Corollary~\ref{coro:creatable}.
		
	\end{enumerate}

\end{definition}

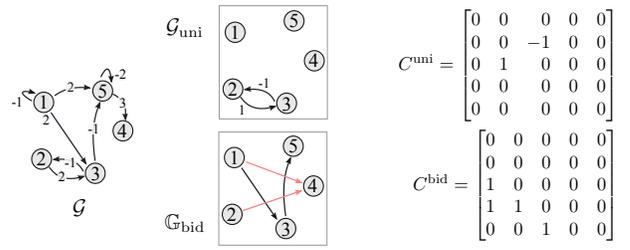
\begin{figure}[t]
 	\centering
 	\begin{tikzpicture}[scale=1]

 		\node[inner sep=0pt, scale=1.05] at (0,0) {\includegraphics{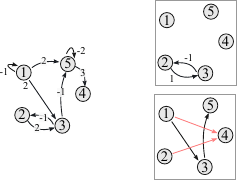}};
 			\node[inner sep=0pt, scale=0.7]  at (4.5,0.8) {$
 			C^{\uni}
 			=
 			\begin{bmatrix*}[r]
 				0 & 0 & 0  & 0  & 0\\
 				0  & 0&  -1& 0  & 0 \\
 				0 & 1 & 0 & 0 & 0 \\ 
 				0  & 0 & 0 & 0  & 0 \\
 				0 & 0 & 0  & 0  & 0 
 			\end{bmatrix*}
 			$};
 			
 			\node[inner sep=0pt, scale=0.7] (f1) at (4.6,-0.8) {$
 				C^{\bid}
 				=
 				\begin{bmatrix*}[r]
 					0 & 0 & 0  & 0  & 0\\
 					0 & 0&  0& 0  & 0 \\
 					1 & 0& 0 & 0 & 0 \\ 
 					1 & 1 & 0 & 0 & 0 \\
 					0 & 0 & 1  & 0 & 0 
                \end{bmatrix*}
 				$};
 				
 			\node[scale=0.8] at (0.2,1.3) {$\CG_{\uni}$};
 			\node[scale=0.8] at (0.2,-1.3) {$\BG_{\bid}$};
 			\node[scale=0.8] at (-1.2,-1.1) {$\CG$};

 	\end{tikzpicture}
 	\caption{Illustration of the unidirectionally and bidirectionally modifiable graphs and their corresponding (un-)weighted adjacency matrices. Red edges indicate the ones that can be created.}
 	\label{modi_map}
 \end{figure}

Note that $\CE_\inc$ and $\CE_\dec$ can have intersections with $\CE_\ctr $. Further, we define a sign graph $\CG_\uni :=(\CV,\CE_\uni,C^\uni )$, where $\CE_\uni=\CE_\inc \cup \CE_\dec$, and $ C^\uni =[c^\uni _{ij}]_{n \times n}$ satisfies 
\begin{align*}
	c^\uni _{ij} =\begin{cases*}
		1, \hspace{8pt}\text{ if } (j,i)\in \CE_\inc,\\
		-1, \text{ if } (j,i)\in \CE_\dec,\\
		0, \hspace{8pt}\text{ otherwise}.
	\end{cases*}
\end{align*}
 In addition, we define an unweighted graph $\BG_\bid =(\CV,\CE_\bid)$ where $\CE_\bid=\CE_\ctr\cup \CE_\cre$. The sign graph $\CG_\uni$ encodes which edges can be modified through direct vibrations. The edge signs indicate the direction of the allowed changes. We refer to  $\CG_\uni$ as \textit{unidirectionally modifiable graph}. The unweighted graph $\BG_\bid$ represents edges that can be adjusted in both directions, which are referred to as \textit{bidirectionally modifiable graph}. We further let $C^\bid =[c^\bid_{ij}]_{n\times n}$ be the unweighted adjacency matrix of  $\BG_\bid $, respectively. Given a digraph $\CG$,  $\CG_\uni$ and $\BG_\bid$ and their corresponding adjacency matrices can be found by following these steps  (see Fig.~\ref{modi_map} for an example).

While $\CG_\uni$ and $\BG_\bid$ identify modifiable edges and the direction of potential functional changes, they do not guarantee the feasibility of achieving arbitrary modifications to any subset of edges within $\CE_\uni \cup \CE_\bid$. This is because realizing specific changes often requires multiple vibrations applied at different network locations, which can induce complex cross-effects as discussed in Section~\ref{controllability}. To formally characterize achievable modifications through vibrational control, we introduce the concept of realizable changes.

\begin{definition}[Realizable perturbation]
	Consider a perturbation matrix $\Delta\in \R^{n \times n}$. We say $\Delta$ is \textit{vibrationally realizable} if there exists a vibrational control matrix $V(t)$ for the system~\eqref{system:vib} such that the weighted adjacency matrix of the functioning network system~\eqref{averaged:main} becomes $\bar A = A+\Delta$. 
\end{definition}

Next, we construct graph-theoretic conditions for realizable perturbations. Let $\CG_\Delta = (\CV,\CE_\Delta,\Delta):=\CG(\Delta)$ be the weighted digraph associated with $\Delta$. We start by considering the simplest scenario where  $\Delta$ has only one non-zero entry. This indicates that we intend to modify a single edge within  $\CG$.  From the previous section, we have the following corollary.

\begin{figure*}[t]
	\centering
	\begin{tikzpicture}[scale=1]		
		\node[inner sep=0pt, scale=1.05] at (0,0) {\includegraphics{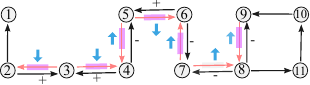}};
		\node[inner sep=0pt, scale=1.05] at (8,0) {\includegraphics{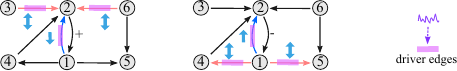}};		
		
		\node[inner sep=0pt, scale=0.7]  at (0,0.9) {(a) Changing multiple directly modifiable edges};	
		
		\node[inner sep=0pt, scale=0.7]  at (7,0.9) {(b) Changing multiple jointly modifiable edges};			
	\end{tikzpicture}
	\caption{Simultaneous modification of multiple edges. (a) Multiple directly modifiable edges: when the target edges do not form a trail longer than 1 (highlighted in red), they can be arbitrarily modified along the directions specified by the bold  blue arrows (increase or decrease) via vibrational inputs. (b) Multiple jointly modifiable edges: when all target edges either leave from or enter at a single node (blue and red edges), and only one of them belongs to a 2-cycle ($(1,2)$ in blue), their weights can be arbitrarily modified along the directions indicated by the bold blue arrows via same-frequency vibrations (a two-way arrow meas an edge can be both increased and decreased). Driver edges are covered by pink blocks.}
	\label{control_multiple}
\end{figure*}

\begin{corollary}[Perturbation on a single edge]\label{con_single}
	Consider a perturbation matrix $\Delta=[\delta_{pq}]_{n\times n}$, where $(j,i)$ is the only edge in $\CG_\Delta$. The perturbation  $\Delta$ is vibrationally realizable if either of the following conditions holds:
	\begin{enumerate}
		\item[(i)]  The edge $(j,i)$ is also an edge in $\CG_\uni$, and its weight satisfies $\sign(\delta_{ij})=c^\uni_{ij}$, 
		\item[(ii)] The edge $(j,i)$ is also an edge in $\BG_\bid$; the opposite edge $(i,j)$ does not exist in $\CG$, i.e., $a_{ji}=0$; and there is a pair of driver edges that do not belong to $\CG_\uni$. 
	\end{enumerate}
\end{corollary}

Note that, the conditions in (ii) prevent the unintended modifications observed in Fig.~\ref{single_edge_stabilization}-(c). Vibrational control inputs in these cases (i) and (ii) can be designed directly following Lemmas~\ref{design:direct} and~\ref{indirect:design}, respectively.

Next, we study the situation where $\Delta$ has multiple  non-zero entries. We first investigate how multiple directly modifiable edges can be changed  simultaneously.


\begin{lemma}[Perturbation on multiple directly modifiable edges]\label{muti_dir}
	Given a perturbation matrix  $\Delta\in \R^{n \times n}$, it is vibrationally realizable if the following conditions are satisfied:
	\begin{enumerate}
		\item[(i)]  The sign graph of $\Delta$, $\CG_\sign(\Delta)=(\CV,\CE_\Delta,S_\Delta)$, satisfies $\CG_\sign(\Delta)\subseteq\CG_\uni$;
		\item[(ii)] The longest trail in $\CG_\sign(\Delta)$ has length 1. 
	\end{enumerate}
	Then, any perturbation matrix $\Delta$ under these conditions can be realized by the vibrational control matrix  $V(t)=[v_{ij}(t)]_{n \times n}$ satisfying \begin{align}\label{design:multi_dir}
		v_{ij}(t)	= \begin{cases}
			u_{ij}\sin (\beta_{ij} t ), \hspace{1.1cm}\forall i,j: \delta_{ij}\neq 0,\\
			0,  \forall t\ge 0, \hspace{1.5cm}\text{ for any other }i,j,
		\end{cases}
	\end{align}
	where the amplitudes and  frequencies are arbitrarily chosen such that $\beta_{ij}$'s are incommensurable and
	$
		\frac{u_{ij}}{\beta_{ij}}=\sqrt{{-2 \delta_{ij}}/{a_{ji}}}
$.
\end{lemma}

Fig.~\ref{control_multiple}-(a) provides an intuitive illustration of the conditions in this lemma. Next, we consider the situation where jointly controllable edges are involved.

\begin{lemma}[Perturbation on multiple jointly modifiable edges]\label{muti_indir}
	Given a perturbation matrix  $\Delta=[\delta_{ij}]\in \R^{n \times n}$, it is vibrationally realizable if the weighted digraph associated with it, $\CG_\Delta=(\CV,\CE_\Delta, \Delta)$, satisfies the following conditions:
	\begin{enumerate}
		\item[(i)] Only one edge in $\CE_\Delta$, denoted $(j_0, i_0)$, belongs to a 2-cycle in $\CG(A)$.  Also, its weight has the same sign as the corresponding edge in $\CG_{\uni}$, i.e., $\sign(\delta_{i_0j_0})=c^\uni_{i_0j_0}$. 		
		\item[(iii)] All the other edges in $\CE_\Delta$ either enter $ i_0$ or leave $j_0$. 
	\end{enumerate}
	
	Then, any perturbation matrix $\Delta$ under these conditions can be realized by the vibrational control matrix  $V(t)=[v_{ij}(t)]_{n \times n}$ satisfying \begin{align}\label{design:multi_bid}
		v_{ij}(t)	= \begin{cases}
			u_{ij}\sin (\beta t ), \hspace{1.1cm}\forall i,j: \delta_{ij}\neq 0,\\
			0,  \forall t\ge 0, \hspace{1.2cm}\text{ for any other }i,j,
		\end{cases}
	\end{align}
	where the amplitudes $u_{ij}$ and the frequency $\beta$ are arbitrarily chosen such that 	\begin{equation}\label{ampli:multi_bid}
		\frac{u_{ij}}{\beta}= \begin{cases*}
			\sqrt{\frac{-2 \delta_{i_0j_0}}{a_{j_0i_0}}}, \hspace{1.5cm} \text{ if } i=i_0, j = j_0 \\
			-\frac{\delta_{ij}}{a_{j_0i_0}}\sqrt{\frac{-2a_{j_0i_0}}{\delta_{i_0j_0}}}, \hspace{0.6cm}\text{for any other } i, j.
		\end{cases*}
	\end{equation}
\end{lemma}

In this lemma, while the direction of change for edge $(j_0, i_0)$ is constrained by its inherent modifiability (increasable or decreasable), the other edges can be modified arbitrarily. Fig.~\ref{control_multiple}-(b) illustrates the constructed conditions.

\begin{figure}[t]
	\centering
	\begin{tikzpicture}[scale=1]		
		\node[inner sep=0pt, scale=1.05] at (0,0) {\includegraphics{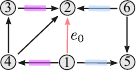}};
		
		\node[inner sep=0pt, scale=0.8]  at (3,0.25) {$\CE^1_\dri =\{(1,4), (3,2)\}$};	
		\node[inner sep=0pt, scale=0.8]  at (3,-0.25) {$\CE^2_\dri =\{(1,5), (6,2)\}$};	
		
	\end{tikzpicture}
	\caption{An illustration of driver sets. To induce a desired modification to the edge $e_0$, vibrational control can be applied to either of the two driver sets,  $\CE^1_\dir$ or $\CE^2_\dri$ (covered by pink blocks and blue blocks, respectively). Both driver sets are different from the target edge $e_0$. }
	\label{multi_driver}
\end{figure} 

\textbf{Driver set:} Given a perturbation matrix $\Delta$, $\CE_\Delta$ denotes the set of target edges to modify. We denote $\CE_\dri(\Delta)$ as the set of edges on which vibrational inputs are imposed to realize the perturbation. We refer to $\CE_\dri(\Delta)$ as the \textit{driver set}. It is worth mentioning that the driver set can be different from the target edges. Moreover, given a perturbation matrix $\Delta$, there can be multiple driver sets that realize it (see Fig.~\ref{multi_driver} for an example).  

\section{Vibrational Stabilizability}\label{stabilizability}

Building upon the previous section's conditions for precise modification of multiple network edges, we now explore the vibrational stabilizability of entire network systems.

\subsection{General Condition}

\begin{theorem}[Vibrational stabilizability]\label{stabilizability_general}
	The system~\eqref{main} is vibrationally stabilizable if there is a perturbation matrix $\Delta\in \R^{n \times n}$ such that 
	\begin{enumerate}
		\item[(i)] the matrix $A+\Delta$ is Hurwitz,
		\item[(ii)] the weighted digraph $\CG_\Delta=(\CV,\CE_\Delta, \Delta)$ can be decompose into $r$ weighted subgraphs, denoted as  $\CG(\Delta^{(p)})=(\CV,\CE_{\Delta^{(p)}}, \Delta^{(p)}), p=1, \dots, r$, such that each $\Delta^{(p)}\in \R^{n\times n}$ satisfies the conditions in Corollary~\ref{con_single}, Lemma~\ref{muti_dir}, or Lemma~\ref{muti_indir},
		\item[(iii)] and there exist  $r$ driver sets, one for each $\CE_{\Delta^{(p)}}$, that are mutually disjoint   
	\end{enumerate}
\end{theorem}

This theorem provides a sufficient condition for vibrational stabilizability of network systems. We provide an algorithm to design vibrational control inputs based on this theorem in Algorithm~\ref{control:placement}. The main idea is to design vibrational inputs for each perturbation $\Delta^{(p)}$ separately, ensuring that the driver sets for different perturbations are mutually disjoint. By selecting incommensurate frequencies for vibrations in different clusters, it is ensured that influences on different clusters are isolated. Consequently, the combined vibrational matrix $V(t)= \sum_{p=1}^{r}V^{(p)}(t)$ achieves the overall desired change $\Delta$. We next demonstrate how Theorem \ref{stabilizability_general} and Algorithm~\ref{control:placement} can be applied by the example below.

\begin{algorithm}[t]
	\caption{Vibrational Control Design}
	\label{control:placement}
	\begin{algorithmic}[1]
		\State \textbf{Input:} System matrix $A$ satisfying Theorem~\ref{stabilizability_general}
		\State Identify a realizable $\Delta$ such that $A+\Delta$ is Hurwitz
		\State Associate $\Delta$ with a weighted digraph $\CG_\Delta$, and decompose $\CG_\Delta$ into $r$ disjoint clusters, $\CG(\Delta^{(p)}), p =1,2,\dots, r$
		\For{$p=1:r$}
		\State	Identify a driver set $\CE_{\rm{dri}}(\Delta^{(p)})$ for each $\CG(\Delta^{(p)})$, ensuring driver sets are mutually disjoint
		\State	Design vibrational inputs following Corollary~\ref{con_single}, Lemma~\ref{muti_dir}, or Lemma~\ref{muti_indir}, while ensuring vibrations introduced to different clusters are incommensurable 
		\EndFor
	\end{algorithmic}
\end{algorithm}

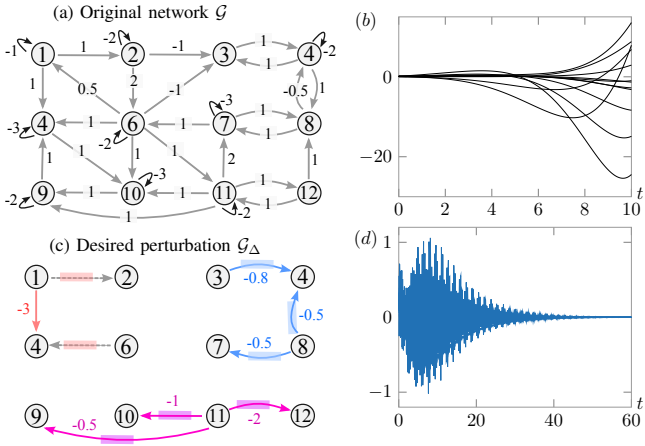
\begin{figure}[t]
	\centering
	\begin{tikzpicture}[scale=1]

		\node[inner sep=0pt, scale=0.7]  at (-0.2,2.8) {(a) Original network $\CG$};
		\node[inner sep=0pt, scale=0.7]  at (2.8,2.7) {$(b)$ };
		
		\node[inner sep=0pt, scale=0.7]  at (0,-0.3) {(c) Desired perturbation $\CG _\Delta$};
		\node[inner sep=0pt, scale=0.7]  at (2.8,-0.2) {$(d)$ };

		\node[inner sep=0pt, scale=1.2] at (0.2,1.3) {\includegraphics{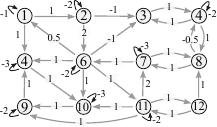}};
		
		\node[inner sep=0pt, scale=1.2] at (0.2,-1.7) {\includegraphics{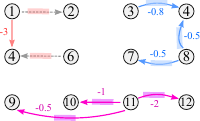}};
		\node[inner sep=0pt, scale=0.7] at (4.6,0) {\includegraphics{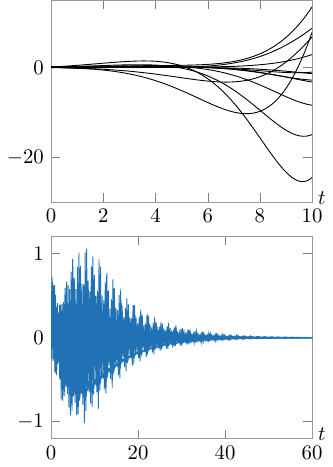}};

	\end{tikzpicture}
	\caption{Illustration of vibrational stabilization. (a) and (b) depict an unstable original network. (c) presents a desired stabilizing perturbation represented by $\CG_\Delta$, decomposed into three clusters (solid edges in distinct colors). Corresponding driver edges (colored blocks) enable the realization of the perturbations in these clusters. (d) demonstrates the stabilizing effect of the vibrational control designed in Example~\ref{vib_general_ctr}.}
	\label{general_stabilizability}
\end{figure}

\begin{example}\label{vib_general_ctr}
	Consider the  system with a weighted digraph depicted in Fig.~\ref{general_stabilizability}-(a). It satisfies the conditions outlined in Theorem~\ref{stabilizability_general}, ensuring the system's vibrational stabilizability.  To stabilize it, we find that it is sufficient to introduce the functional changes described  by the weighted digraph $\CG_\Delta$ in Fig.~\ref{general_stabilizability}-(c).  To enable that, we choose three sets of driver edges that are mutually disjoint. Then, we design the vibrational control inputs using the results in Section~\ref{multiple_contr}. We use joint vibrations $v_{21}(t) =\sin(t)$ and $v_{56}(t) = 3\sin(t)$  to functionally change the weight of $(1,5)$. To functionally modify $(3,4), (8,4),$ and $(8,7)$, we let $v_{43} (t)= \sqrt{3.2}\sin(\sqrt{2}t), v_{48}(t)= \sqrt{3}\sin(\sqrt{3}t),$ and $v_{78}(t)= \sqrt{15}\sin(\sqrt{5}t)$.  To change $(11,9), (11,10),$ and $(11,12)$, we let $v_{12,11} (t)= 2\sqrt{7}\sin(\sqrt{7}t), v_{10,11}(t) = \sqrt{7}\sin(\sqrt{7}t), 	v_{9,11}(t) = 0.5\sqrt{7}\sin(\sqrt{7}t)$. Further, let $\varepsilon=0.04$. As shown in Fig.~\ref{general_stabilizability}-(d), these control inputs  stabilize the network system. 
\end{example}

We remark that, determining whether a system satisfies the conditions in Theorem~\ref{stabilizability_general} may not be always straightforward as it requires identifying a realizable perturbation to stabilize the network. To address this, we investigate structural  stabilizability and construct a more easily verifiable condition by introducing an additional assumption in the next subsection.

\subsection{Structural Condition}

\begin{assumption}\label{stable:scalor}
	Assume that each individual subsystem in~\eqref{main} is stable, i.e., $a_{ii}<0$ for all $i=1,2,\dots, n$.
\end{assumption}

Interconnections between subsystems can result in instability of the overall system, despite individual stability. To stabilize such a system, our strategy is to functionally remove edges from the network by imposing vibrational control. Let us define  removable edge sets.

\begin{definition}[Removable edge set]\label{defi:remove_graph}
An edge set $\CE_1\subseteq \CE$ is said to be \emph{removable} if there exists a vibrational control such that the unweighted graph of the functioning system~\eqref{averaged:main} becomes $\BG(\bar A) = (\CV, \CE \setminus \CE_1)$. 
\end{definition}

Note that, unlike the precise modification of edge weights without affecting other edges, functionally removing edges allows for changes to the weights of other edges as long as no new edges are created. Denote the sets of directly, path-enabled, and 2-cycle-enabled removable edges in $\CG$ as $\CE_\rem^{\dir}, \CE_\rem^{\pat}$, and $\CE_\rem^{\cyc}$, respectively.
Leveraging the results in Section~\ref{multiple_contr}, we have the next corollary for removable edge sets. 

\begin{corollary}\label{Coro:removable}
 The edge set $\CE_1\subseteq \CE$ is \emph{vibrationally removable} if one of the following situations are satisfied:
	\begin{enumerate}
		\item[(i)] $\CE_1$ only contains a single edge $(j,i)$, which satisfies any of the conditions: (a) it is directly removable, i.e., $(j,i) \in \CE_\rem^{\dir}$;
		(b) it is path-enabled removable but not directly modifiable (i.e., $(j,i) \in \CE_\rem^{\pat}$ and  $(j,i)\notin \CE_\uni$);
		(c) it is 2-cycle enabled removable (e.g., $(j,i) \in \CE_\rem^{\cyc}$).
		\item[(ii)] $\CE_1\subseteq \CE_\rem^{\dir}$, and its longest trail has length 1.
		\item[(iii)] There is an edge $(j_0, i_0) \in\CE_1$ such that:  $(j_0, i_0)$ is the only edge that belongs to a 2-cycle;  $(j_0, i_0)\in \CE_\rem^\dir$; the other edges either all enter the node $i_0$ or all leave $j_0$.
		\item[(iv)] None of the edges in $\CE_1$ belongs to a 2-cycle, and there exists an edge $(j_0, i_0)\in\CE/\CE_1$ that belongs to a 2-cycle such that the edges in  $\CE_1$  either all enter the node $i_0$ or all leave $j_0$.
	\end{enumerate}	
\end{corollary}

\begin{theorem}[Structural vibrational stabilizability]\label{structural}
Let $\BG=(\CV,\CE)$ be the unweighted digraph of the system~\eqref{main:compact}. The system~\eqref{main} is vibrationally stabilizable if there are $r$ weakly connected subgraphs in $\BG$, denoted as $\BG^{(p)}=(\CV,\CE^{(p)}), p =1, 2, \dots, r$, such that 
\begin{enumerate}
	\item[(i)] Each edge set $\CE^{(p)}$ satisfies one of the four situations in Corollary \ref{Coro:removable},
	\item[(ii)] Removing the edges in these subgraphs from $\BG$ results in a directed acyclic graph\footnote{A DAG is a directed graph that does not contain any directed cycles.} (DAG),
	\item[(iii)] There exist $r$  driver sets,  one for removing each $\CE^{(p)}$, that are mutually disjoint.
\end{enumerate}
\end{theorem}

	 Theorem~\ref{structural} indicates that  a network system can be stabilized by functionally removing edges so that the remaining network is a DAG.  The main intuition is to   eliminate feedback loops that could introduce destabilizing effects.
Vibrational inputs can be designed following Algorithm~\ref{control:placement}, but focusing on removing edges in $\BG$. We next provide an example to illustrate how this theorem can be applied.

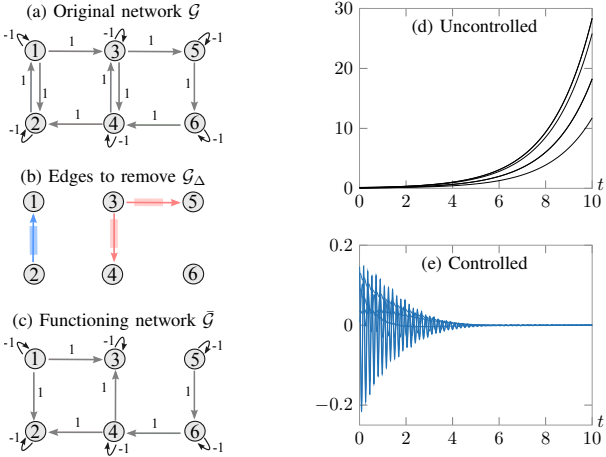
\begin{figure}[t]
	\centering
	\begin{tikzpicture}[scale=1]

		\node[inner sep=0pt, scale=0.7]  at (-0.6,1.9) {(a) Original network $\CG$};
		\node[inner sep=0pt, scale=0.7]  at (-0.6,-.3) {(b) Edges to remove $\CG_\Delta$};
		\node[inner sep=0pt, scale=0.7]  at (-0.6,-2.2) {(c) Functioning network $\bar \CG$};
		\node[inner sep=0pt, scale=1.05] at (-0.6,0.9) {\includegraphics{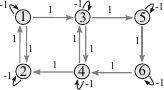}};
		
		\node[inner sep=0pt, scale=1.05] at (-0.6,-1.1) {\includegraphics{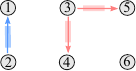}};
		
		\node[inner sep=0pt, scale=1.05] at (-0.6,-3.2) {\includegraphics{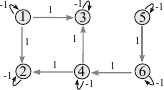}};

		\node[inner sep=0pt, scale=0.7] at (4,-0.9) {\includegraphics{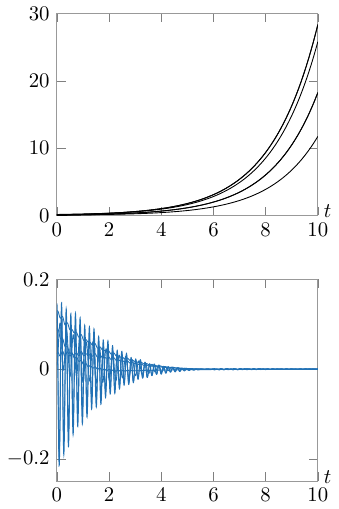}};
		
		\node[inner sep=0pt, scale=0.7]  at (4.2,1.7) {(d) Uncontrolled};
		
		\node[inner sep=0pt, scale=0.7]  at (4.2,-1.45) {(e) Controlled};
	\end{tikzpicture}
	\caption{Illustration of structural stabilizability. (a) Original unstable network. (b) Edges to remove and the corresponding driver sets (covered by colored blocks, forming into clusters). (c) Removing the edges in (b) results in a directed acyclic graph. The originally unstable system (d) becomes stable under vibrational inputs to the driver edges (e).}
	\label{structural_stabilizable}
\end{figure}

\begin{example}\label{example:structural}
	Consider a linear system with a weighted digraph depicted in Fig.~\ref{structural_stabilizable}-(a). Despite individual subsystems being stable, the overall system is not (shown in Fig.~\ref{structural_stabilizable}-(d)).   Following  Theorem~\ref{structural}, one can identify two weakly connected subgraphs shown in Fig.~\ref{structural_stabilizable}-(b). Removing these edges results in a DAG (see Fig.~\ref{structural_stabilizable}-(c)). In addition, to functionally remove them, the driver edge sets $\CE_\dir^{(1)}=\{(2,1)\}$, and $\CE_\dir^{(2)}=\{(3,4), (3,5)\}$  are disjoint, indicating the system is vibrationally stabilizable. 
	
	 We next follow similar steps in Algorithm~\ref{control:placement} to design  the vibrational control matrix $V(t)=[v_{ij}(t)]$  in~\eqref{system:vib}. Particularly, we use sinusoidal inputs, where $v_{ij}(t)=u_{ij}\sin(\beta_{ij}t)$. As the edge $(2,1)$ is directly removable, following Lemma~\ref{design:direct} we select $u_{12} =\sqrt{2}$ and $\beta_{12} = 1$. To functionally remove the edges $(3,4)$ and $(3,5)$ jointly, we utilize Lemma~\ref{muti_indir} to design vibrational inputs. Specifically, we let	$u_{43}=u_{53}= 2, \beta_{43} = \beta_{53}= \sqrt{2}$. Also, let $\varepsilon=0.04$. It is worth noting that same-frequency vibrations have been introduced to the edges $(3,4)$ and $(3,5)$. As shown in Fig.~\ref{structural_stabilizable}-(e), these vibrational inputs stabilize the system, demonstrating our findings in Theorem~\ref{structural}.
\begin{remark}
	When designing vibrational inputs, selecting the value of $\varepsilon$ is also important. Lemma \ref{lemma:averaged} proves the existence of a threshold, below which any $\varepsilon$ suffices. Similar to many results built on averaging theory (e.g., see~\cite[Th. 10.4]{HKK:02-bis}), determining the exact value of the threshold is difficult. In Examples \ref{vib_general_ctr} and \ref{example:structural}, we find that a small value of $\varepsilon$ on the order of $0.01$ is sufficient. In general, a careful selection of $\varepsilon$, tailored to the system under consideration, is required. 
\end{remark}

\end{example}

\section{Conclusion}

Vibrational control offers a distinct advantage by eliminating the need for sensors to measure a system's internal states. This work applies this open-loop strategy to stabilize complex networks. We demonstrate that vibrational control operates by modifying edge weights on average within network systems. We have established graph-theoretic conditions to identify modifiable edges and determine overall network stabilizability. Additionally, we present methods for designing vibrational inputs to achieve desired network changes.

Our findings open avenues for future research. Notably, the presented conditions for vibrational stabilizability are sufficient but not necessary, necessitating further investigation to close the gap. Current vibrational control methods rely on fixed amplitude and frequency inputs. Incorporating a feedback loop to dynamically adjust these parameters could optimize system regulation. Moreover, the current vibrations introduce only a $1/\varepsilon$ term in the system, requiring only first-order averaging techniques. It is worth exploring vibrations that generate higher-order terms $1/\varepsilon^n$ for $n>1$. Analyzing such systems requires higher-order averaging techniques (see~\cite{SAV:01,SJA-FMJ:07}), which may lead to interesting insights. In addition, our results are limited to systems where each node dynamics is one-dimensional. Extending our results to higher-dimensional dynamics requires alternative analytical tools, presenting an intriguing direction for future research. Furthermore, given the potential connections to brain stimulation therapies, our theoretical framework could inform the development of enhanced and predictable treatments for brain disorders through synergistic collaborations with experimental research.

\appendix

\subsection{Analysis of Section of \ref{joint_control}}
We only prove Lemma~\ref{indirect:design}, since  Corollary~\ref{coro:creatable} can be shown by following similar steps. 

\begin{pfof}{Lemma~\ref{indirect:design}}
	By assumption in this lemma, there exist a directed trail $\{(j,p), (p,q), (q,i)\}$. We construct the proof by exhausting all the three situations in~\eqref{sit_1}-\eqref{sit_3}.
	
	We first show the case (i), where the nodes $j,p,q, i$ are all distinct. Without loss of generality, we  reorder the nodes in the network such that $j = 1$, $p=2$, $q=3$, and $i=4$. Now, it remains to show the edge $(1,4)$ can be functionally modified. The vibrational inputs in~\eqref{theo:amp} are then introduced to the edges $(1,2)$ and $(3,4)$, leading to the vibrational control matrix 	
	\begin{equation}
		V(t)=\scalemath{0.8}{\begin{bmatrix}
			0&0&0& 0& 0 &0\\
			v_{21}(t)&0&0& 0& 0&0\\
			0&0&0& 0& 0&0\\
			0 &0&v_{43}(t)& 0& 0&0\\
			\vdots&\vdots&\vdots&\vdots& \ddots&\vdots\\
			0&0&0& 0& 0&0\\
		\end{bmatrix}},
	\end{equation}
	where $v_{21}(t)=u_{21}\sin(\beta t)$ and $v_{43}(t)=u_{43}\sin(\beta t)$. Now, consider an  auxiliary system  described by
	$
		d{\hat x}/ds=V(s)\hat x.
	$
	Then, its fundamental matrix $\Psi(s)$ can be derived as
	\begin{equation}
		\Psi(s)=\scalemath{0.8}{\begin{bmatrix}
			1&0&0& 0& 0&0\\
			\psi_{21}(s)&1&0& 0& 0&0\\
			0&0&1& 0& 0&0\\
			0 &0&\psi_{43}(s)& 1& 0&0\\
			\vdots&\vdots&\vdots&\vdots&\ddots &\vdots\\
			0&0&0& 0& 0&1\\
		\end{bmatrix}},
	\end{equation}
	with $\psi_{21}(s)=-u_{21}/\beta\cdot \cos(\beta s)$ and $\psi_{43}(s)=-u_{43}/\beta \cdot \cos(\beta s)$. 
 As in Section~\ref{functioning}, following $\bar A= \lim_{T\to \infty} \frac{1}{T}\int_{t=0}^{T} \Psi^{-1}(t)A\Psi(t) \dd t$, one can derive that  $\bar A=[\bar a_{k\ell}]$ of the averaged system~\eqref{averaged:main} satisfies 
	\begin{align*}
		\bar a _ {41} &= a_{41}+ \lim_{T\to \infty} \frac{1}{T}\int_{t=0}^{T} \psi_{43}(t) a_{31} -\psi_{43}(t)a_{32}\psi_{21}(t) \dd t\\
		& = a_{41} -\frac{1}{2\beta^2}  u_{43} a_{32} u_{21} . 
	\end{align*}
	Given that $a_{32}\neq 0$, to induce a functional change $\delta\in \R$ to the edge $(1,4)$, one simply needs to select $u_{43}, u_{21}$, and $\beta$ such that $\delta = -\frac{1}{2\beta^2}  u_{43} a_{32} u_{21}$.
	
	Next, we prove the case (ii), where the directed trail is $\{(j,i), (i,q), (q,i)\} $. Without loss of generality, we reorder the nodes in the network such that $j=1, i =3,$ and $q=2$. The vibrational inputs in~\eqref{theo:amp} are then introduced to the edges $(1,3)$ and $(2,3)$, leading to the vibrational control matrix 	
	\begin{equation}
		V(t)=\scalemath{0.8}{\begin{bmatrix}
			0&0&0& 0& 0 \\
			0&0&0& 0& 0\\
			v_{31}(t)&v_{32}(t)&0& 0& 0\\
			\vdots&\vdots&\vdots&\ddots& \vdots\\
			0&0&0& 0& 0\\
		\end{bmatrix}},
	\end{equation}
	with $\psi_{31}(s)=-{u_{31}}/{\beta} \cdot \cos(\beta s)$ and $\psi_{32}(s)=-{u_{32}}/\beta  \cdot \cos(\beta s)$.	
	Following similar steps as above, one can derive that $\bar A=[\bar a_{k\ell}]$ of the averaged system satisfies
	$
		\bar a_{31} = a_{31} - \frac{1}{2\beta^2} u_{32} a_{23} u_{31}   - \frac{u_{31}^2}{2\beta^2} a_{13}.
	$
	For any $\delta \in \R$, selecting $u_{31}, u_{32}$, and $\beta$ such that 
	$
		- \frac{1}{2\beta^2} u_{32} a_{23} u_{31}   - \frac{u_{31}^2}{2\beta^2} a_{13}=\delta
	$
	leads to a functional change to the edge $(1,3)$ by $\delta$. 
	
	Finally, we show the case (iii), where the directed trail is $\{(j,p), (p,j), (j,i)\}$. We introduce vibrations to the edges $(j,p)$ and $(j,i)$. Following similar steps as above, one can show that vibrational inputs $v_{pj}(t)	=u_{pj}\sin (\beta t )$ and $v_{ij}(t)	= u_{ij}\sin (\beta t )$ lead to an averaged system~\eqref{averaged:main} with $\bar A = [\bar a_{k\ell}]$ satisfying 
	$
		\bar a_{ij} = a_{ij} - \frac{1}{2\beta^2} u_{ij} a_{jp} u_{pj}   - \frac{u_{ij}^2}{2\beta^2} a_{ji}
	$. 
	Then, for any $\delta\in\R$, one can always select $u_{ij}$, $u_{pj} $, and $\beta$ such that 
	$
		- \frac{1}{2\beta^2} u_{ij} a_{jp} u_{pj}   - \frac{u_{ij}^2}{2\beta^2} a_{ji}=\delta
	$
	a functional change $\delta$ is induced to the edge $(j,i)$. The proof is complete.
	\end{pfof}
	
\subsection{Proof of Lemma~\ref{muti_dir}}
Since the longest trail in $\CG_\Delta$ has length $1$, each node is either a sink,  source or isolated. Assume that there are $m_1$ source nodes, $m_2$ sink nodes, and $m_3=n-m_1-m_2$ isolated nodes. There exists a permutation matrix $P$ such that 
\begin{equation*}
	\Delta'=P\Delta P^{-1} = \begin{bmatrix}
		\Bfzr_{m_1\times m_1}&\Bfzr_{m_1\times m_2} & \Bfzr\\
		\hat \Delta & \Bfzr_{m_2\times m_2} &\Bfzr\\
		\Bfzr&\Bfzr&\Bfzr
	\end{bmatrix},
\end{equation*}
where $\hat \Delta=[\hat \delta_{ij}]\in \R^{m_2\times m_1}$. Let $A'=PA P^{-1}$, and rewrite $A'$ into a block-matrix form
$
	A':= \begin{bmatrix}
		*&\tilde{A} &*\\
		\hat A & *&*\\
		*&*&*
	\end{bmatrix},
$
where $\hat A =[\hat a_{ij}] \in \R^{m_2\times m_1}$ and $\tilde A=[\tilde a_{ij}] \in \R^{m_1\times m_2}$.

To show that there exists a vibration matrix such that $\bar A=A +\Delta$, it suffices to show that there exists $V'(t)$ such that the averaged system of the controlled system 
$\dot x = \left(A'+\frac{1}{\varepsilon}V'\left(\frac{t}{\varepsilon}\right)\right)x$, denoted as $\dot {\bar x} = {\bar A}' \bar x$, satisfies $\bar A=A' +\Delta'$.

To this end, consider that the vibrational control matrix $V'(t)$ has the block-matrix form
\begin{equation*}
	V'(t) = \begin{bmatrix}
		\Bfzr_{m_1\times m_1}&\Bfzr & \Bfzr\\
		\hat V(t) & \Bfzr_{m_2\times m_2} &\Bfzr\\
		\Bfzr&\Bfzr&\Bfzr
	\end{bmatrix},
\end{equation*}
with $\hat V(t)=[\hat v_{ij}(t)]\in \R^{m_2\times m_1}$ satisfying 
$
	\hat v_{ij}(t)=  \beta_{ij}\sqrt{\frac{-2 \delta_{ij}}{\tilde a_{ji}}} \sin (\beta_{ij}t),\text{ if } \hat \delta_{ij}\neq 0,
$
and $\hat v_{ij}(t)=0$ for all $t\ge 0$, otherwise. Following similar steps as in Section~\ref{functioning}, one can derive the fundamental matrix of the $dx/ds = V'(s)x$ is 
\begin{equation*}
	\Psi'(s) = \begin{bmatrix}
		I_{m_1\times m_1}&\Bfzr & \Bfzr\\
		\hat \Psi(s) & I_{m_2\times m_2} &\Bfzr\\
		\Bfzr&\Bfzr&I
	\end{bmatrix},
\end{equation*}
where $\hat \Psi(s)=[\hat \psi_{ij}]$ satisfies 
$
	\hat \psi_{ij}(s)=  - \sqrt{\frac{-2 \hat \delta_{ij}}{\tilde a_{ji}}} \cos (\beta_{ij}s)$, if $ \hat  \delta_{ij}\neq 0,
$
and $\hat \psi_{ij}(s)=0$ otherwise. Note that Condition (i) implies that $-\hat \delta_{ij}/\tilde a_{ji}>0$, ensuring that $\sqrt{-\hat \delta_{ij}/\tilde a_{ji}}$ is well defined. Then, it can be derived that 
\begin{equation*}
	{\Psi'}^{-1}(s) = \begin{bmatrix}
		I_{m_1\times m_1}&\Bfzr & \Bfzr\\
		-\hat \psi(s) & I_{m_2\times m_2} &\Bfzr\\
		\Bfzr&\Bfzr&I
	\end{bmatrix}.
\end{equation*}
Subsequently, from Section~\ref{functioning}, it holds that $\bar A'= \lim_{T\to \infty} \frac{1}{T}\int_{t=0}^{T} {\Psi'}^{-1}(t)A'\Psi'(t)\dd t$, which can be derived as 
\begin{align*}
	\bar A'= A'+ \begin{bmatrix}
		\Bfzr_{m_1\times m_1}&\Bfzr&\Bfzr\\
		 H & \Bfzr_{m_2\times m_2}&\Bfzr\\
		\Bfzr&\Bfzr&\Bfzr
	\end{bmatrix},
\end{align*}
where $ H=[h_{ij}]\in\R^{m_2\times m_1}$ satisfies
$
	 h_{ij}=\lim\limits_{s\to \infty}\frac{1}{T}\int_{0}^{T}- \frac{2 \hat \delta_{ij}}{\tilde a_{ji}} \cos^2 (\beta_{ij}s)\cdot \tilde a_{ji} ds=\hat \delta_{ij}.
$
This implies that $H= \hat\Delta$, which completes the proof.

\subsection{Proof of Lemma~\ref{muti_indir}}

We construct the proof by consider the two situations separately: (i) the edges in $\CE_\Delta$ all enter the node $i_0$, and (ii) the edges in $\CE_\Delta$ all leave the node $j_0$. Denote $e_0\coloneq (j_0,i_0)$.

We start with the situation (i).  Assume that the number of edges in $\CE_\Delta/{e_0}$ is $m$. Then, there exists a permutation matrix $P$ such that 
\begin{equation*}
	\Delta'=P\Delta P^{-1} = \begin{bmatrix}
		\Bfzr_{m-1}&\Bfzr& \cdots&\Bfzr&\Bfzr&\Bfzr\\
		\delta'_{m1} & \delta'_{m2} &\cdots &\delta'_{m,m-1} &0&\Bfzr\\
		\Bfzr_{n-m}&\Bfzr&\Bfzr&\Bfzr&\Bfzr&\Bfzr
	\end{bmatrix},
\end{equation*}
where $\delta'_{m1}$ corresponds to the desired change to the directly modifiable edge $(j_0,i_0)$.  Applying the same permutation to the matrix $A$ results in 
\begin{equation*}
	A'=PA P^{-1} = \begin{bmatrix}
			*&a'_{1m}&*\\
			* & \Bfzr_{m-2}&*\\
			*&*&*
		\end{bmatrix}.
\end{equation*}
Next, we show that there exists $V'(t)$ such that the averaged system of the controlled system 
$\dot x = \left(A'+\frac{1}{\varepsilon}V'\left(\frac{t}{\varepsilon}\right)\right)x$, denoted by $\dot {\bar x} = {\bar A}' \bar x$, satisfies $\bar A=A' +\Delta'$.

To this end, we consider the  vibrational control matrix
\begin{equation*}
	V' (t)= \begin{bmatrix}
		\Bfzr_{m-1}&\Bfzr& \cdots&\Bfzr&\Bfzr&\Bfzr\\
		v'_{m1}(t) & v'_{m2}(t)  &\cdots &v'_{m,m-1}(t)  &0&\Bfzr\\
		\Bfzr_{n-m}&\Bfzr&\Bfzr&\Bfzr&\Bfzr&\Bfzr
	\end{bmatrix},
\end{equation*}
where 	$v'_{mk}(t)=u_{mk} \sin (\beta t), k=1,2, \dots, m-1$ satisfy
\begin{align*}
	\frac{u_{m1}}{\beta}&=   \sqrt{\frac{-2 \delta'_{m1}}{a'_{1m}}} \sin (\beta t),\\
	\frac{u_{mk}}{\beta}&=  -\frac{\delta'_{mk}}{a'_{m1}} \sqrt{ \frac{-2a'_{m1}}{a'_{1m}} }\sin (\beta t), \text{ for }k=2,\dots, m-1.
\end{align*}
Then, the fundamental matrix of the system $dx/ds=(A'+V'(s))x$ is 
\begin{equation*}
	\Psi'(t) = \begin{bmatrix}
		&&I_{m-1}&&\Bfzr&\Bfzr\\
		\psi'_{m1}(t) & \psi'_{m2}(t)  &\cdots &\psi'_{m,m-1}(t)  &1&\Bfzr\\
		\Bfzr_{n-m}&\Bfzr&\Bfzr&\Bfzr&\Bfzr&I
	\end{bmatrix},
\end{equation*}
where \begin{align*}
	\psi'_{m1}(t)&=   \sqrt{\frac{-2 \delta'_{1m}}{a'_{m1}}} \cos (\beta t),\\
	\psi'_{mk}(t)&=  -\frac{\delta'_{mk}}{a'_{m1}} \sqrt{ \frac{-2a'_{m1}}{a'_{1m}} }\cos (\beta t).
\end{align*}
Then, from $\bar A'= \lim_{T\to \infty} \frac{1}{T}\int_{t=0}^{T} {\Psi'}^{-1}(t)A'\Psi'(t)\dd t$  we have 
\begin{align*}
	\bar a'_{mk} = a'_{mk}-\lim\limits_{T\to \infty}\frac{1}{T}\int_{t=0}^{T}\psi'_{mk}(t) a'_{1m} \psi'_{m1}(t)=a'_{mk}+\delta'_{mk}
\end{align*}
for all $k=1, 2, \dots, k$, implying that $\bar A'=A'+\Delta'$. 

Following similar steps as above, one can show the case (ii) and complete the proof.

\subsection{Proof of Theorem~\ref{stabilizability_general}}

By assumption, we have $\Delta= \Delta^{(1)} + \Delta^{(2)}+\dots+\Delta^{(r)}$ and the driver sets for them are mutually disjoint, then there exists a permutation matrix $P$ such that 
$
	\Delta'=P\Delta P^{-1} =\blk\left(	\Delta'_1, 	\Delta'_2,\dots,	\Delta'_r,\Bfzr\right),
$
where each $\Delta'_i$ has zero entries at diagonal and upper-diagonal positions and corresponds to the perturbations defined by $\Delta^{(i)}$. Since each $\Delta^{(i)}$ satisfies the conditions in  Corollary~\ref{con_single}, Lemma~\ref{muti_dir}, or Lemma~\ref{muti_indir},  there exists a vibrational control matrix of the block-diagonal form $\bar V^{'(i)} = A'+\blk(\Bfzr,\dots,V'_i(t),\dots,\Bfzr)$ such that the averaged matrix $\dot {\bar x}=\bar A^{'(i)} \bar x$ of the controlled matrix $dx'/ds = (A'+V^{'(i)}(s))x$ satisfies 
$
	\bar A^{'(i)} = A'+\blk\left(\Bfzr,\dots,\Delta'_i,\dots,\Bfzr\right).
$
Now, consider the vibrational control matrix 
\begin{equation*}
	V'(t)= \blk\left(V'_1(t),V'_2(t),\dots,V'_r(t),\Bfzr\right),
\end{equation*}
where $V'_i(t)$'s are chosen such that they mutually have incommensurable frequencies.
The fundamental matrix associated with $V'(t)$ is 
\begin{equation*}
	\Psi'(t)= \blk\left(\Psi'_1(t),\Psi'_2(t),\dots,\Psi'_r(t),I\right).
\end{equation*}
Using the fact that the $\Psi'_i(t)$'s have incommensurable frequencies, one can derive that 
$$
	\bar A'= \lim_{T\to \infty} \frac{1}{T}\int_{t=0}^{T} {\Psi'}^{-1}(t)A'\Psi'(t)\dd t=A'+ \Delta',
$$
which completes the proof.

\subsection{Proof of Theorem~\ref{structural}}

\begin{lemma}\label{lemma:stable_graph}
	Assume that the network system described in~\eqref{main} associated with the unweighted   graph $\BG=(\mathcal{V}, \mathcal{E})$. Then, this system is asymptotically stable if the graph $\mathcal{G}$ is a DAG.
\end{lemma}

	\begin{pfof}{Lemma~\ref{lemma:stable_graph}}
	As $\CG$ is also a DAG, according to~\cite{JBJ-GG:00},  it can be topologically ordered. Therefore, one can  arrange the nodes of $\CG$ as a linear ordering that is consistent with all edge directions. In other words, there exists a permutation matrix $P$ such that the matrix $\hat A=:PAP^{-1}$ is lower-triangular. Since $a_{ii}<0$ for all $i$, one can derive that the diagonal entries of $\hat A$ are all negative. This means that $\hat A$ is Hurwitz, implying that $A$ is also Hurwitz. 	
\end{pfof}

To prove Theorem~\ref{structural}, one can follow similar steps as those for Theorem~\ref{stabilizability_general}. One simply needs to construct vibrations to remove the corresponding edges instead of modifying them.

\section{Reference}
\bibliographystyle{IEEEtran}
\bibliography{alias,refers}

\end{document}